# Construction of Benford Random Variables: Generators and Seed Functions

Frank A. Benford

## 1. Introduction

In this paper I describe two general methods of constructing a probability density function (pdf) of a Benford random variable, and show that the pdf of every Benford random variable may be constructed by these methods. These constructions are helpful in exploring the nature of Benford random variables, and may be used to investigate some issues of base dependency of Benford random variables.

A cursory familiarity with Benford random variables might lead one to think that the pdf of a base 10 Benford random variable necessarily has the form

$$f(x) = \frac{\Lambda}{x} \qquad \text{for all} \quad x \in [a, a \times 10^r] \tag{1.1}$$

and $f(x) = 0$ outside of this interval, where $\Lambda$ is a suitable constant, $a$ is a positive real number, and $r$ is a positive integer, or that $f$ is pieced together with parts of this form defined over disjoint domains. This assumption is incorrect. In fact, the pdf of Benford random variables may have a large variety of forms, including forms that are positive and continuously differentiable over all of $\mathbb{R}_{++}$, the set of strictly positive real numbers, and have the form given in eq. (1.1) over no interval of $\mathbb{R}_{++}$ of positive length.

The form given in eq. (1.1) has some geometric "clues" that it's the pdf of a base 10 Benford random variable. One of my aims when I began this research was to find ways to construct pdfs of a base $b$ Benford random variable that have no obvious geometric features that reveal $b$.

## 2. Benford Random Variables.

We begin with a review of some basic facts about Benford random variables. For a comprehensive summary of the mathematical theory of Benford random variables, the reader is referred to the recent text by Berger and Hill [3].

Historically, the "Benford property" of random variables (and functions and sequences) was defined in terms of the joint distribution of the most significant digits of numbers. More recently, it has been found convenient to define the Benford property in terms of the *significand* function. Berger and Hill define the significand of $x$ for any $x \in \mathbb{R}$, but for our purposes it's sufficient to restrict the definition to numbers $x > 0$.

Let $b > 1$. **Definition 2.1**. For any $x > 0$, the *base* $b$ *significand of* $x$, denoted $S_b(x)$, is defined as the unique $s \in [1, b)$ such that

$$x = s \times b^k$$

for some (necessarily unique) $k \in \mathbb{Z}$. Hence

$$x = S_b(x) \times b^k \quad \text{where} \quad S_b(x) \in [1, b) \text{ and } k \in \mathbb{Z}. \tag{2.1}$$

For example, $2020 = 2.02 \times 10^3$, so $S_{10}(2020) = 2.02$. The reader may show that $S_8(2020) = 3.9453125$.

Now let $X$ be a positive random variable; that is, $\Pr(X > 0) = 1$. Assume that $X$ is continuous with a probability density function (pdf). **Definition 2.2**: $X$ is base $b$ Benford (or $X$ is $b$-Benford) if and only if the cumulative distribution function (cdf) of $S_b(X)$ is given by

$$\Pr(S_b(X) \leq s) = \log_b(s) \quad \text{for all} \quad s \in [1, b). \tag{2.2}$$

It's useful at this point to introduce some non-standard notation. Let $y \in \mathbb{R}$ and recall that the "floor" of $y$, written $\lfloor y \rfloor$, is defined as the largest *integer* that is less than or equal to $y$. Define $\langle y \rangle$ as:

$$\langle y \rangle \equiv y - \lfloor y \rfloor \tag{2.3}$$

and note that $0 \leq \langle y \rangle < 1$ for every $y \in \mathbb{R}$. We'll call $\langle y \rangle$ the *fractional part* of $y$, though if $y < 0$ this description is misleading.

If we take the logarithm base $b$ of eq. (2.1) we obtain

$$\log_b(x) = \log_b(S_b(x)) + k. \tag{2.4}$$

On the other hand,

$$\log_b(x) = \langle \log_b(x) \rangle + \lfloor \log_b(x) \rfloor. \tag{2.5}$$

As $\lfloor \log_b(x) \rfloor$ is necessarily an integer and $0 \leq \log_b(S_b(x)) < 1$, comparison of eqs. (2.4) and (2.5) shows that

$$\log_b(S_b(x)) = \langle \log_b(x) \rangle \qquad \text{and} \qquad k = \lfloor \log_b(x) \rfloor \tag{2.6}$$

for any $x > 0$.

Using eq. (2.6), we may rephrase Definition 2.2 in several logically equivalent ways: $X$ is $b$-Benford iff

(1) $\Pr(\log_b(S_b(X)) \leq \log_b(s)) = \log_b(s)$ for all $s \in [1, b)$,
(2) $\Pr(\langle \log_b(X) \rangle \leq u) = u$ for every $u \in [0, 1)$,
(3) $\langle \log_b(X) \rangle \sim U[0, 1)$,
(4) $X = b^Y$ where $\langle Y \rangle \sim U[0, 1)$,

where the notation "$W \sim U[0, 1)$" means that $W$ is uniformly distributed on the half open interval $[0, 1)$. (More generally, I use the symbol " $\sim$ " to mean "is distributed as." Hence, for example, "$X \sim f()$" means that $X$ is distributed with pdf $f()$, "$X_1 \sim X_2$" means that $X_1$ and $X_2$ have the same distribution, and "$X \sim N(\mu, \sigma^2)$" means that $X$ is distributed normally with mean $\mu$ and variance $\sigma^2$.)

For any random variable $Y$, if $\langle Y \rangle \sim U[0,1)$ we sometimes say that $Y$ is "uniformly distributed modulo one," abbreviated "u.d. mod 1." Hence $X$ is $b$-Benford iff $\log_b(X)$ is u.d. mod 1.

To find algorithms that allow one to create $b$-Benford random variables, we first consider the converse problem of determining whether a given random variable is or is not $b$-Benford. Given $b > 1$ and the positive random variable $X$, define $Y \equiv \log_b(X)$ so $X = b^Y$. Let

$$\begin{array}{ll} f() & \text{denote the pdf of } X, \\ g() & \text{denote the pdf of } Y, \text{ and} \\ \widetilde{g}() & \text{denote the pdf of } \langle Y \rangle. \end{array}$$

Then $X$ is $b$-Benford iff $\widetilde{g}(u) = 1$ for almost all $u \in [0,1)$.

A standard argument shows that $g()$ may be obtained from $f()$ and vice-versa by the formulas

$$g(y) = \ln(b) b^y f(b^y) \qquad \text{and} \qquad f(x) = \frac{\Lambda_b}{x} g(\log_b(x)) \tag{2.7}$$

where $\Lambda_b \equiv 1/\ln(b)$. Under a reasonable assumption[1], we may show that

$$\widetilde{g}(u) = \sum_{k \in \mathbb{Z}} g(k+u). \tag{2.8}$$

**Example 2.1**. Let $k \in \mathbb{Z}$ be fixed, and suppose that $X$ is distributed on $\mathbb{R}_{++}$ with the pdf given by

$$f(x) = \begin{cases} \Lambda_b / x & \text{if} \quad b^k \le x < b^{k+1}, \\ 0 & \text{elsewhere.} \end{cases}$$

Then eq. (2.7) yields

$$g(y) = \begin{cases} 1 & \text{if} \quad k \le y < k+1, \\ 0 & \text{elsewhere.} \end{cases}$$

As $\widetilde{g}(u) = g(k+u) = 1$ for any $u \in [0,1)$, it follows that $X$ is $b$-Benford.

**Example 2.2**. If $f()$ is a convex combination of pdfs of the form given in Example 2.1 and $X \sim f()$, then $X$ is $b$-Benford. To state this formally, let $(\alpha_k)_{k \in \mathbb{Z}}$ be a sequence of non-negative numbers such that

[1] Explained in some detail in [2]. The series on the right-hand side of eq. (2.8) is viewed as a series of functions mapping $u \in [0,1)$ into $\mathbb{R}$, and we assume that the convergence of this series is uniform.

$$\sum_{k\in\mathbb{Z}} \alpha_k = 1. \tag{2.9}$$

Note that the collection of intervals $([b^k, b^{k+1}))_{k\in\mathbb{Z}}$ partitions $\mathbb{R}_{++}$, and define $f\colon \mathbb{R}_{++} \to \mathbb{R}_+$ as

$$f(x) = \frac{\alpha_k \Lambda_b}{x} \qquad \text{for } x \in [b^k, b^{k+1}) \tag{2.10}$$

for all $k \in \mathbb{Z}$. Then

$$g(y) = \alpha_k \qquad \text{for } y \in [k, k+1) \tag{2.11}$$

for all $k \in \mathbb{Z}$. Hence, for all $u \in [0,1)$,

$$\breve{g}(u) = \sum_{k\in\mathbb{Z}} g(k+u) = \sum_{k\in\mathbb{Z}} \alpha_k = 1,$$

so $X$ is $b$-Benford whenever $X \sim f()$.

**Example 2.3 (Berger and Hill, p. 32, Example 3.6 (iii)).** Suppose that $X$ is distributed on $\mathbb{R}_{++}$ with the following pdf:

$$f(x) = \begin{cases} \Lambda_{10}\left(\frac{1}{x} - \frac{1}{x^2}\right) & \text{if} \quad x \in [1, 10), \\ 10\Lambda_{10} \cdot \frac{1}{x^2} & \text{if} \quad x \in [10, 100), \\ 0 & \text{elsewhere.} \end{cases} \tag{2.12}$$

Applying eq. (2.7) with $b = 10$ to this pdf yields

$$g(y) = \begin{cases} 1 - 10^{-y} & \text{if} \quad 0 \le y < 1, \\ 10^{1-y} & \text{if} \quad 1 \le y < 2, \\ 0 & \text{elsewhere.} \end{cases}$$

Here's a picture of $g()$.

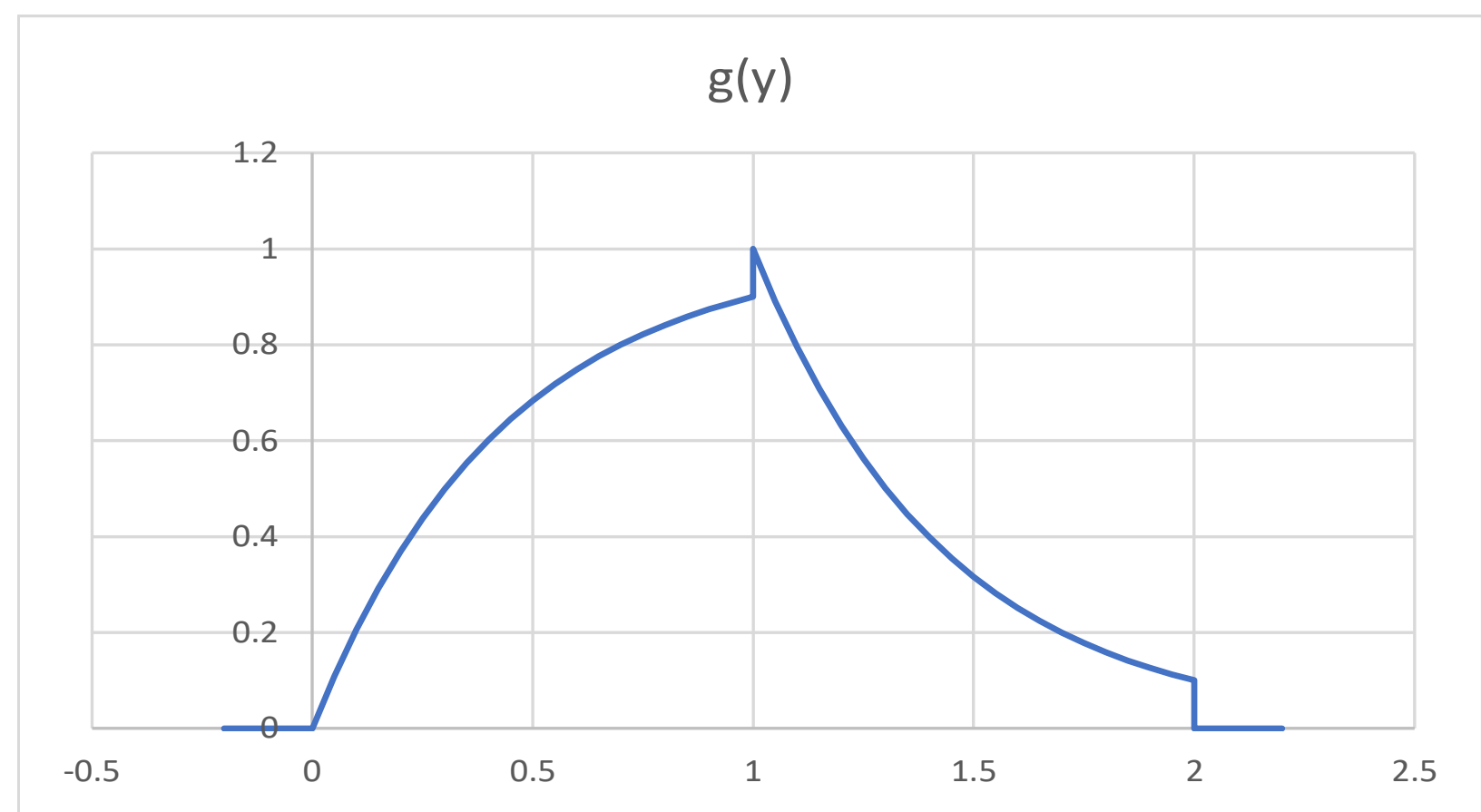

As

$$\sum_{k=0}^{1} g(k+u) = (1 - 10^{-u}) + 10^{-u} = 1,$$

it follows that $X$ is 10-Benford.

Note that there are discontinuities in this graph at $y = 1$ and $y = 2$, and that the support of $Y$ is precisely the interval $[0, 2]$. These geometric features of $g()$ are inherited by the associated graph of $f()$ and provide clues to the base $b$.

## 3. Generators and the Construction of Benford Random Variables.

In this section I describe a general algorithm for the creation of the pdf of Benford random variables using a vector of functions called a "generator." In fact, the pdf of every Benford random variable may be generated by this algorithm.

**Definition 3.1**. A vector of functions $(g_k())_{k\in\mathbb{Z}}$ each mapping $[0, 1)$ into $\mathbb{R}_+$ is called a **generator**. A generator is said to be **uniform** if

$$\sum_{k\in\mathbb{Z}} g_k(u) = 1 \tag{3.1}$$

for all $u \in [0, 1)$, where we assume, as above, that the convergence of the series in eq. (3.1) is uniform.

The only generators of interest to us in this document are those associated with a pdf $g()$ by the formula

$$g_k(u) = g(k + u) \tag{3.2}$$

for every $k \in \mathbb{Z}$. This pdf may be retrieved from the associated generator by the formula

$$g(y) = g_{\lfloor y \rfloor}(\langle y \rangle) \tag{3.3}$$

for every $y \in \mathbb{R}$. If a random variable $Y \sim g()$ and the generator associated with $g()$ is uniform, then $\langle Y \rangle \sim U[0, 1)$ so $b^Y$ is $b$-Benford.

Every uniform generator may be thought of as partitioning the unit square into "strata," where the height of the $k$-th stratum is given by $g_k(u)$.

**Example 3.1**. The partition of the unit square for Example 2.3 is as shown below. Here $g_0(u) = 1 - 10^{-u}$, $g_1(u) = 10^{-u}$, and $g_k(u) = 0$ for every other $k$.

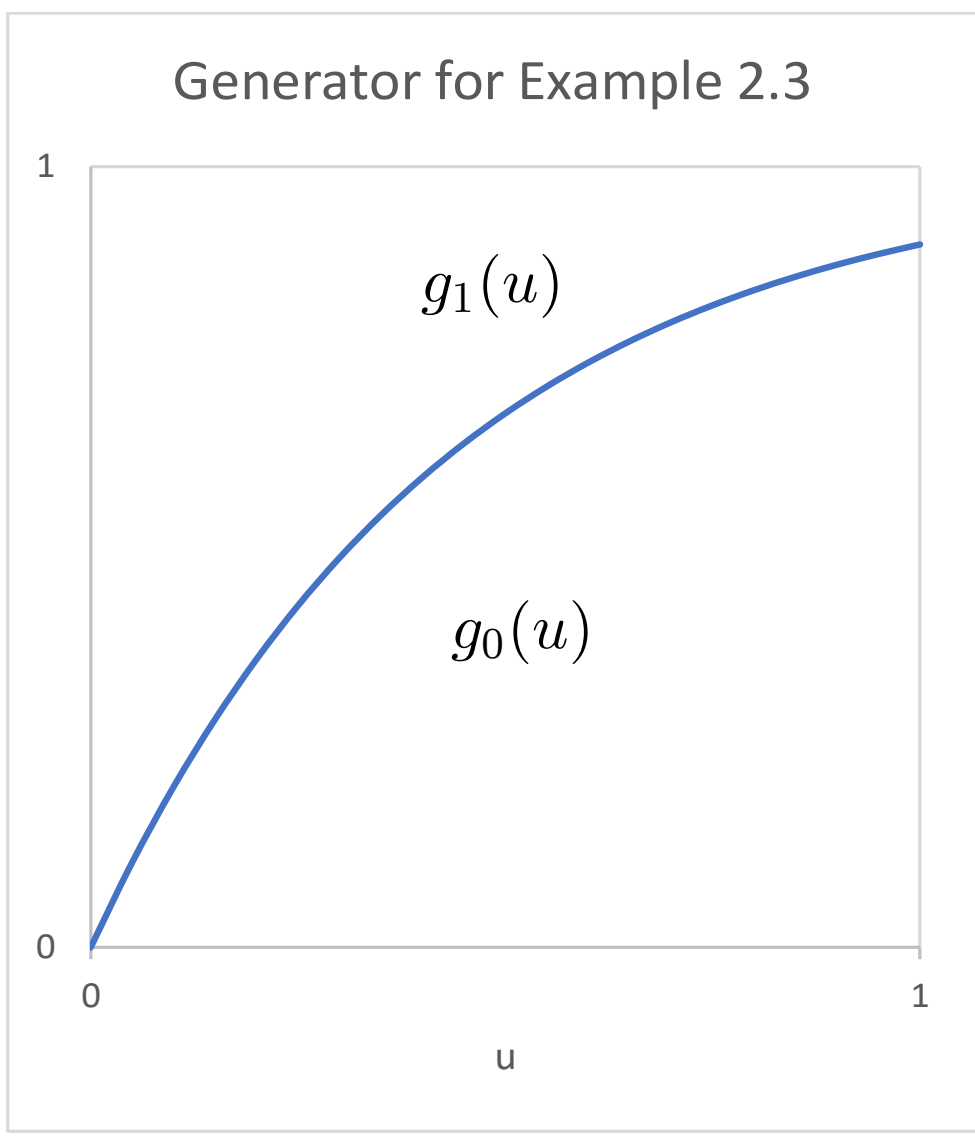

**Example 3.2**. Consider the generator with $g_k(u) \equiv 0$ for all $k$ except $k = 0, 1,$ and 2, and

$$
\begin{aligned}
g_0(u) &\equiv \sin^2(\alpha u), \\
g_1(u) &\equiv \sin^2[\alpha(1+u)] - \sin^2(\alpha u), \\
g_2(u) &\equiv 1 - \sin^2[\alpha(1+u)]
\end{aligned}
$$

where $\alpha \equiv \pi/4$. This generator is clearly uniform, and the associated partition of the unit square is shown in the following figure.

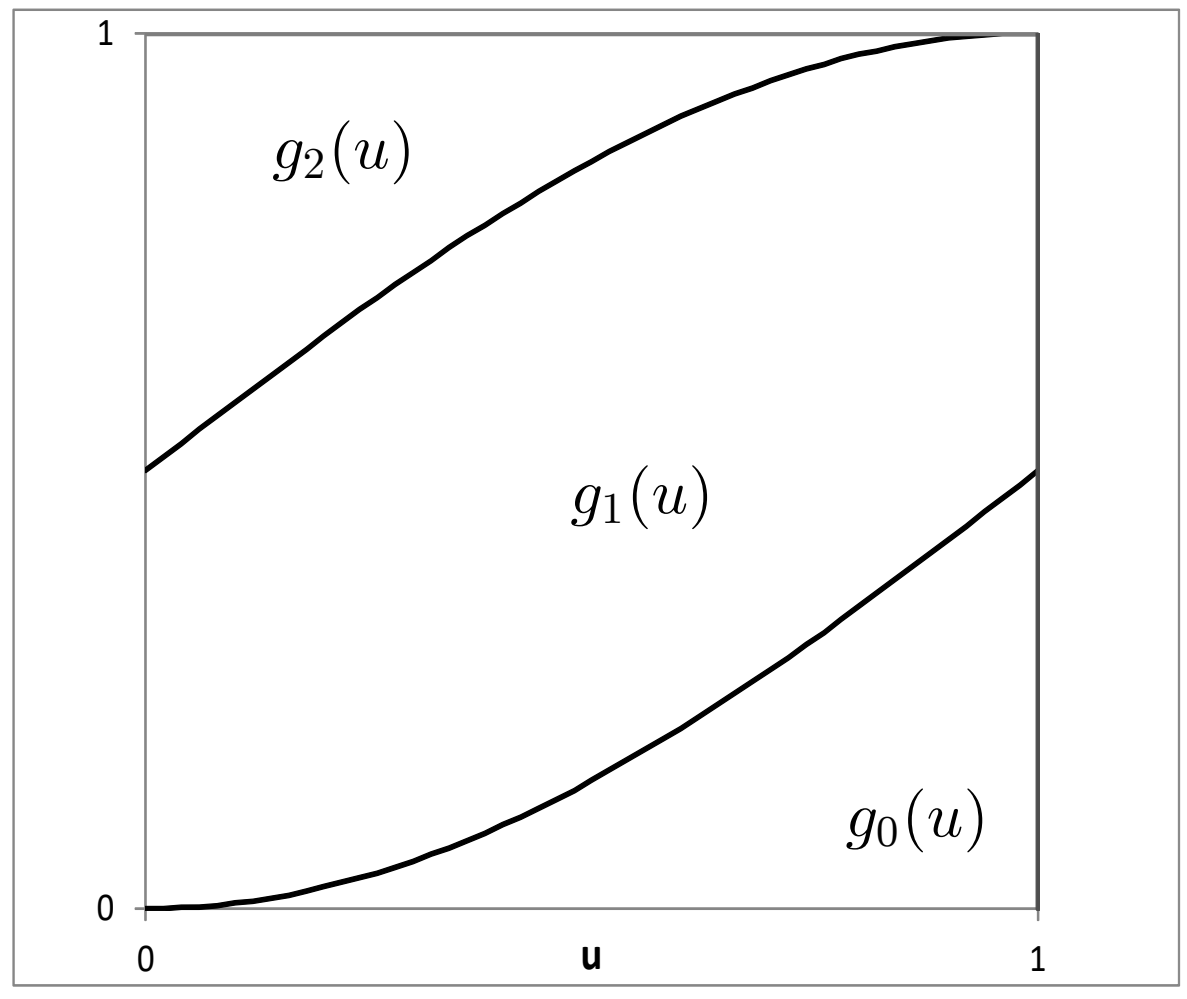

## 4. Seed Functions.

There are an infinite number of ways to produce a uniform generator and its associated pdf $g()$. In this section I describe a method based on what I call a "seed function."

**Definition 4.1**. A **seed function** is a function $H\colon \mathbb{R} \to [0,1]$ with the following three properties: (1) $H()$ is "unit interval increasing," by which I mean

$$H(y-1) \leq H(y) \qquad \text{for all} \quad y \in \mathbb{R}, \tag{4.1}$$

(2) $H(y) \to 0$ as $y \to -\infty$, and (3) $H(y) \to 1$ as $y \to +\infty$.

If $H()$ is increasing[2] then it is unit interval increasing, but not conversely in general, so every cdf is a seed function, but not every seed function is a cdf. A seed function may or may not be continuous everywhere on its domain, and it may or may not be differentiable on some or all of the points where it's continuous.

Let $H()$ be a seed function and define $g\colon \mathbb{R} \to \mathbb{R}_+$ as

$$g(y) \equiv H(y) - H(y-1) \tag{4.2}$$

for every $y \in \mathbb{R}$. The unit interval increasing condition implies that $g(y) \geq 0$ for all $y \in \mathbb{R}$. The graph of $H(y-1)$ is just the graph of $H(y)$ shifted to the right by 1. It may be shown that

$$\int_{-\infty}^{\infty} g(y)\, dy = 1. \tag{4.3}$$

This is geometrically obvious if $H()$ is an increasing function: fill up the space between $H(y)$ and $H(y-1)$ with horizontal line segments, each of length 1. Hence, $g()$ is a legitimate probability density function.

**Proposition 4.1**: The generator associated with $g()$ is uniform. **Proof**. Fix $u \in [0,1)$ and consider the sum

$$\sum_{k=-N}^{N} g(k+u) = \sum_{k=-N}^{N} [H(k+u) - H(k-1+u)] = H(N+u) - H(-N-1+u)$$

for some integer $N \geq 1$. Letting $N \to \infty$, we find

$$\sum_{k \in \mathbb{Z}} g(k+u) = 1,$$

which completes the proof.

**Corollary 4.2**. If $g()$ is given by eq. (4.2) where $H()$ is a seed function, and if $g()$ is the pdf of a random variable $Y$, then $Y$ is u.d. mod 1 and $X \equiv b^Y$ is $b$-Benford.

The generator associated with this $g()$ is constructed by the following rule: for every $k \in \mathbb{Z}$ and $u \in [0,1)$,

[2] Throughout this paper I use "increasing" in a loose sense, as a synonym for "non-decreasing." Hence, I'll say (for example) that $H()$ is *increasing* iff $H(y_1) \leq H(y_2)$ whenever $y_1 < y_2$. If $H(y_1) < H(y_2)$ whenever $y_1 < y_2$, I'll say that $H()$ is *strictly* increasing. A similar convention applies to "decreasing."

$$g_k(u) \equiv g(k+u) = H(k+u) - H(k-1+u).$$

What we've done, intuitively, is wrap the seed function around a cylinder of circumference 1. The number $u$ measures the fraction of the circumference traveled in the current cycle, and the function $g_k(u)$ measures the increase in $H()$ from the same position in the previous cycle.

The construction of $g()$ from $H()$ is reversible, and $H()$ may be retrieved from $g()$ as follows: for any $y \in \mathbb{R}$,

$$\begin{aligned} H(y) &= g(y) + H(y-1) \\ &= g(y) + g(y-1) + H(y-2) \\ &\vdots \\ &= \sum_{\ell=0}^{\infty} g(y-\ell). \end{aligned} \tag{4.4}$$

Hence there is a one-to-one map between seed functions and pdfs with uniform generators.

**Example 4.1.** The generator used in Example 3.2 was derived from the seed function

$$H(y) \equiv \begin{cases} 0 & \text{if} \quad y < 0, \\ \sin^2((\pi/4)y) & \text{if} \quad 0 \le y < 2, \\ 1 & \text{if} \quad y \ge 2. \end{cases} \tag{4.5}$$

**Example 4.2**. The seed function

$$H(y) = \begin{cases} 0 & \text{if} \quad y < 0, \\ 1 - 10^{-y} & \text{if} \quad 0 \le y < 1, \\ 1 & \text{if} \quad y \ge 1 \end{cases}$$

generates the pdf of Example 2.3.

One of my concerns when I began these investigations was to find pdfs for $b$-Benford random variables that offered no geometric clues about $b$. The sort of clues I had in mind include discontinuities, "kinks" (i.e., discontinuities in the derivative), and bounds on the support. Suppose that $X$ is $b$-Benford with pdf $f()$. Then $g()$, the pdf of $\log_b(X)$, is related to $f()$ through eq. (2.7), and there exists a seed function $H()$ associated with $g()$ through eqs. (4.2) and (4.4). The geometric properties of $f()$ are inherited from $H()$ through $g()$. Hence, if $H()$ is continuously differentiable and strictly increasing on $\mathbb{R}$, then $f()$ is continuously differentiable on $\mathbb{R}_{++}$ and $f(x) > 0$ for all $x \in \mathbb{R}_{++}$.

**Example 4.3**. Suppose $H(y) = \Phi(y)$, the cdf of a standard normal random variable. Then $g(y) \equiv H(y) - H(y-1)$ looks like this:

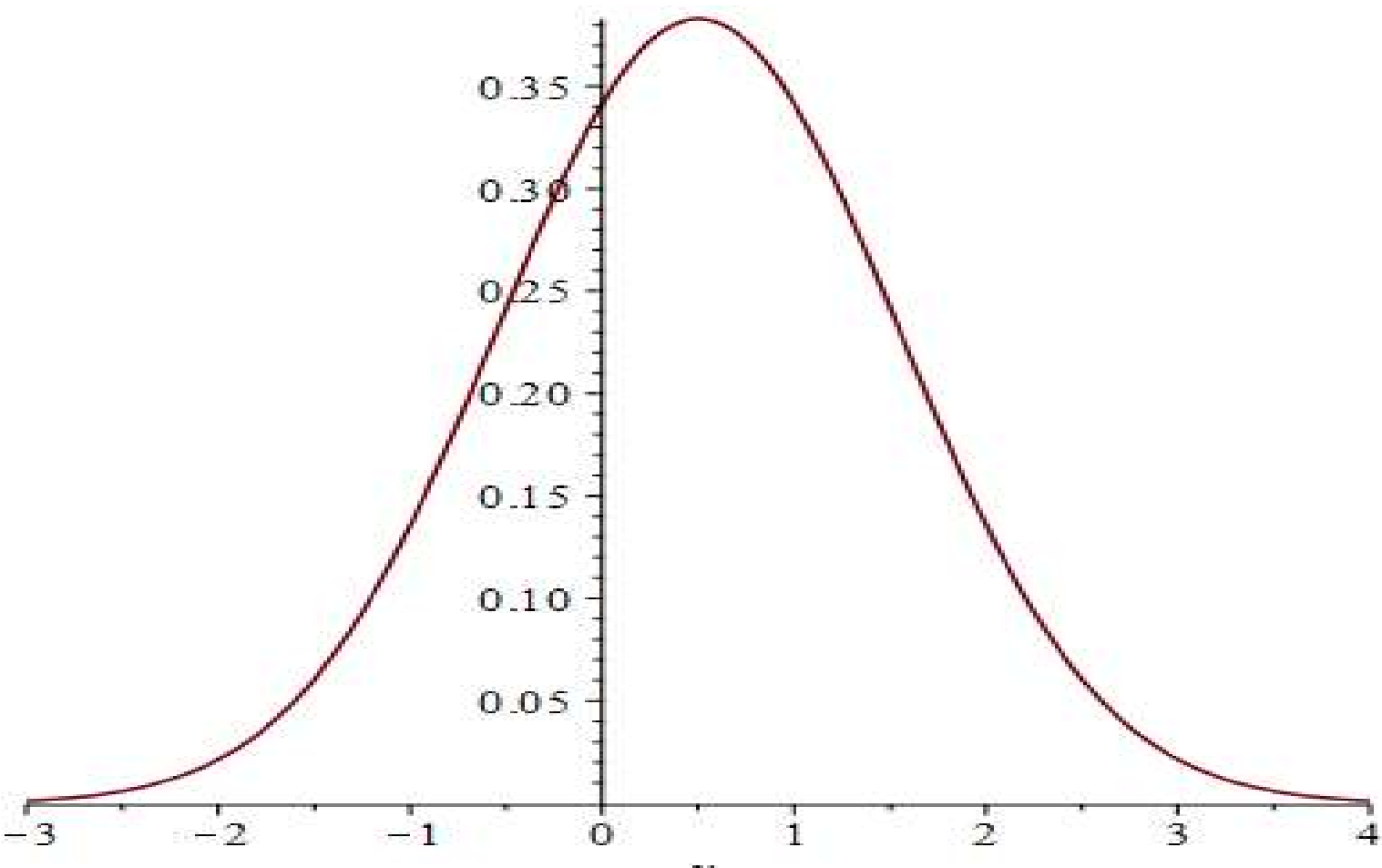

That is, $g()$ is closely approximated by the pdf of a normal random variable with mean $1/2$ and variance 13/12. In fact, the difference between $g()$ and a normal pdf with this mean and variance never exceeds $0.0004$ in absolute value. The implied graph of $f()$ is visually indistinguishable from the pdf of a lognormal random variable. As $g()$ is continuously differentiable and positive for all $y \in \mathbb{R}$, it follows that $f()$ is continuously differentiable and positive for all $x \in \mathbb{R}_{++}$.

The relationships between $f()$, $g()$, $(g_k())_{k\in\mathbb{Z}}$, and $\tilde{g}()$ may be illustrated as follows:

$$f() \overset{b}{\leftrightarrow} g() \leftrightarrow (g_k())_{k\in\mathbb{Z}} \to \tilde{g}()$$

The doubled headed arrows in this diagram indicate that a construction is reversible. The only non-reversible construction in this diagram is from $(g_k())_{k\in\mathbb{Z}}$ to $\tilde{g}()$. If $\tilde{g}(u) = 1$ for almost all $u \in [0, 1)$, the random variable $X$ with pdf $f()$ is Benford. The only place the base $b$ appears in this figure is in the constructions from $f()$ to $g()$ and back. If $g()$ is constructed from a seed function $H()$, the associated generator is always uniform and the diagram is modified as follows:

$$\begin{array}{c} H() \\ \updownarrow \\ f() \overset{b}{\leftrightarrow} g() \leftrightarrow (g_k())_{k\in\mathbb{Z}} \to \tilde{g}() \equiv 1 \end{array}$$

## 5. Shifts and Scalings.

In this section we restrict our attention to increasing seed functions.

**Proposition 5.1**. Suppose that $H_0()$ is an increasing seed function and $\phi\colon \mathbb{R} \to \mathbb{R}$ satisfies the following conditions:

(1) $\phi()$ is increasing,
(2) $\phi(y) \to -\infty$ as $y \to -\infty$, and
(3) $\phi(y) \to +\infty$ as $y \to +\infty$.

Then $H_1(y) \equiv H_0(\phi(y))$ is an increasing seed function. The proof is left to the reader.

For example, we might have $\phi(y) \equiv y^3$ or $\phi(y) = \sinh(y) = \frac{1}{2}(e^y - e^{-y})$. **In this section, however, we'll restrict our attention to the simple "scale and shift" operation $\boldsymbol{\phi(y) \equiv py - q}$ with $\boldsymbol{p > 0}$.** That is, we'll investigate the seed function $H_1(y) \equiv H_0(py - q)$. Let $g_0()$ and $f_0()$ be associated with $H_0()$ through eqs. (4.2) and (2.7) with some base $b$, let $g_1()$ and $f_1()$ be similarly associated with $H_1()$, let $Y_0 \sim g_0()$, $X_0 \sim f_0()$, $Y_1 \sim g_1()$, and $X_1 \sim f_1()$. How are $g_1()$, $Y_1$, $f_1()$, and $X_1$ related to $g_0()$, $Y_0$, $f_0()$, and $X_0$ under this scale and shift operation? We'll restrict our attention to two special cases.

**Case 1: $\boldsymbol{p = 1}$.** In this case, $H_1(y) \equiv H_0(y - q)$, so

$$\begin{aligned} g_1(y) &= H_1(y) - H_1(y-1) \\ &= H_0(y-q) - H_0(y-q-1) = g_0(y-q). \end{aligned}$$

That is, the pdf $g_1()$ is just the pdf $g_0()$ shifted to the right by $q$, and hence $Y_1 = Y_0 + q$. For any $b > 1$, $X_1 \equiv b^{Y_1} = b^q b^{Y_0} = cX_0$ is $b$-Benford, where $c \equiv b^q$. In summary, a shift in the argument of a seed function leads to a shift in the argument of the associated pdf $g()$, and hence to a rescaling of the associated $b$-Benford random variable $X$. Beginning with the scale factor $c$, we may retrieve the associated shift by $q = \log_b(c)$.

This result is equivalent to the following well known fact (e.g., Berger and Hill [3], Theorem 6.3): let $X$ be a positive random variable and let $c > 0$. If $X$ is $b$-Benford, then $cX$ is also $b$-Benford.

**Case 2: $\boldsymbol{q = 0,\ p \neq 1}$.** In this case, $H_1(y) = H_0(py)$, so

$$g_1(y) = H_1(y) - H_1(y-1) = H_0(py) - H_0(py - p). \tag{5.1}$$

This pdf is not functionally related to $g_0(py) = H_0(py) - H_0(py-1)$, so scaling the argument of $H_0()$ does not lead to scaling the argument of $g_0()$. However, the function $g_0(py)$ is related to $X_0$ and $f_0()$. Suppose that $Y_0 \sim g_0()$ and $X_0 \equiv b^{Y_0}$, so $Y_0 = \log_b(X_0)$. We know that $X_0$ is $b$-Benford. Suppose we now wish to examine the Benford properties of $X_0$ relative to another base, say $a > 1$. Let

$$Y \equiv \log_a(X_0) = \frac{\ln(X_0)}{\ln(a)} = \frac{\ln(b)}{\ln(a)} \cdot \frac{\ln(X_0)}{\ln(b)} = p^{-1} Y_0 \tag{5.2}$$

where $p \equiv \ln(a)/\ln(b)$. Let $g()$ denote the pdf of $Y$. The reader may show that

$$g(y) = p g_0(py). \tag{5.3}$$

**Note**: although the random variable $Y_0 \sim g_0()$ is u.d. mod 1 by construction, it is not generally true that $Y$ defined by eq. (5.2) is u.d. mod 1. That is, a random variable $X$ that is

$b$-Benford is not generally $a$-Benford when $a \neq b$. **Example**. Suppose that $Y_0 \sim U[0,1)$, so

$$g_0(y) = \begin{cases} 1 & \text{if} \quad y \in [0,1), \\ 0 & \text{elsewhere.} \end{cases}$$

Also, suppose that $a = b^2$, so $p = 2$. Then

$$g(y) = 2g_0(2y) = \begin{cases} 2 & \text{if} \quad y \in [0,1/2), \\ 0 & \text{elsewhere.} \end{cases}$$

Hence, $Y$ is not u.d. mod 1.

## 6. Examples of Seed Functions.

In this section I present some simple seed functions and examine the implied pdfs $g()$ and $f()$.

**(1) Step functions**. The simplest possible seed function is a simple step function:

$$H(y) = \begin{cases} 0 & \text{if} \quad y < q, \\ 1 & \text{if} \quad y \geq q, \end{cases}$$

for some fixed $q \in \mathbb{R}$. This seed function yields a square pdf $g()$:

$$g(y) = \begin{cases} 1 & \text{if} \quad q \leq y < q+1, \\ 0 & \text{elsewhere.} \end{cases}$$

This seed function jumps from 0 to 1 in a single step. More generally, we may consider step functions that increase from 0 to 1 in multiple steps. Suppose we are given a countable set of "jump points" $\{q_k\}$, a set of "indicator" step functions linked to the jump points

$$I(y \geq q_k) \equiv \begin{cases} 0 & \text{if} \quad y < q_k, \\ 1 & \text{if} \quad y \geq q_k, \end{cases}$$

and a countable set of "jump sizes" $\{\alpha_k\}$ such that $\alpha_k > 0$ for all $k$ and $\sum_k \alpha_k = 1$. We may then define $H()$ as

$$H(y) = \sum_k \alpha_k I(y \geq q_k).$$

The jump points need not be separated by integral amounts. For example, consider the seed function $H(y)$ with jump points $(0, 0.75, 1.5)$ and jump sizes $(0.2, 0.5, 0.3)$ shown below:

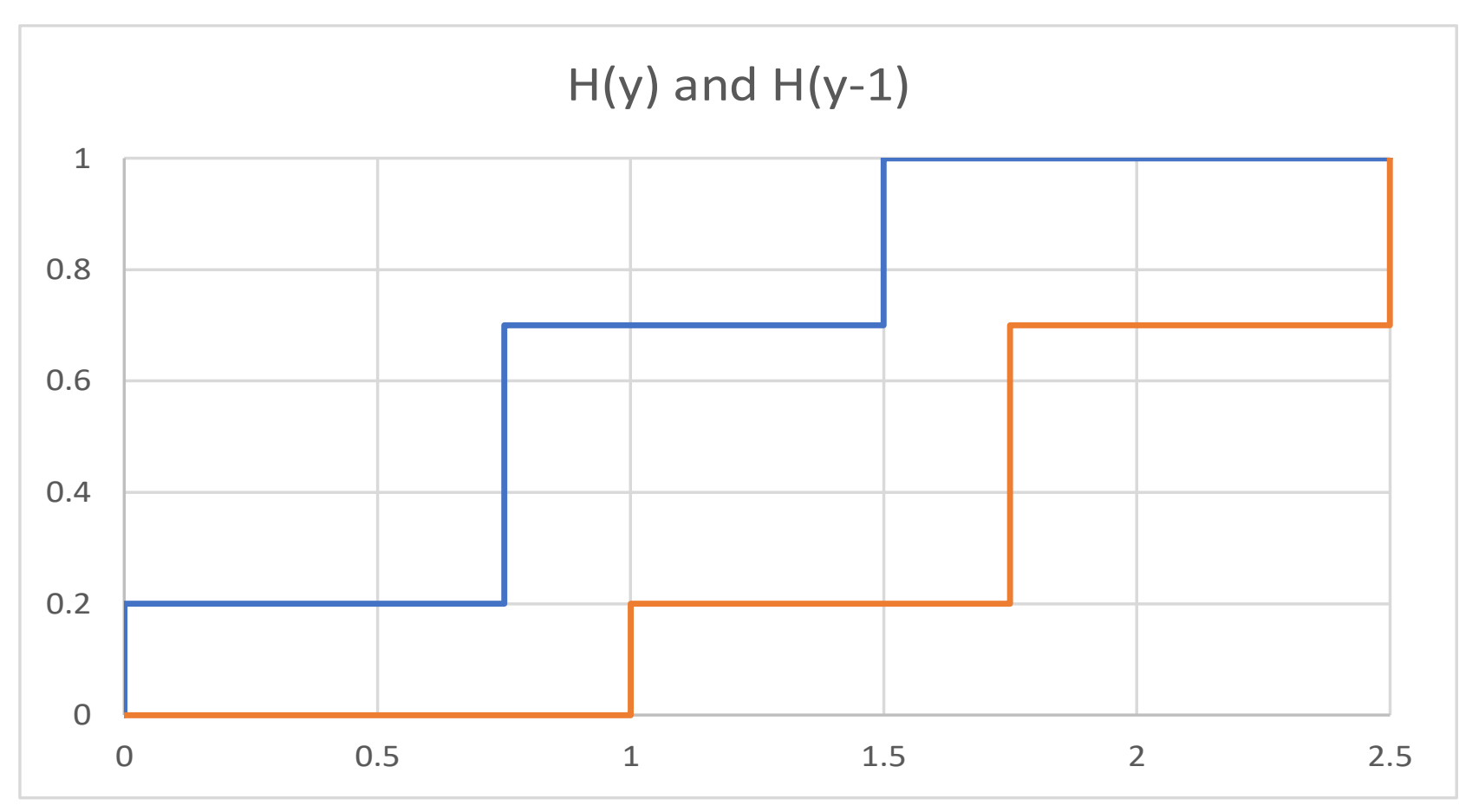

The “shifted” seed function $H(y-1)$ is shown in orange in this figure. This seed function implies that $g()$ is as shown below

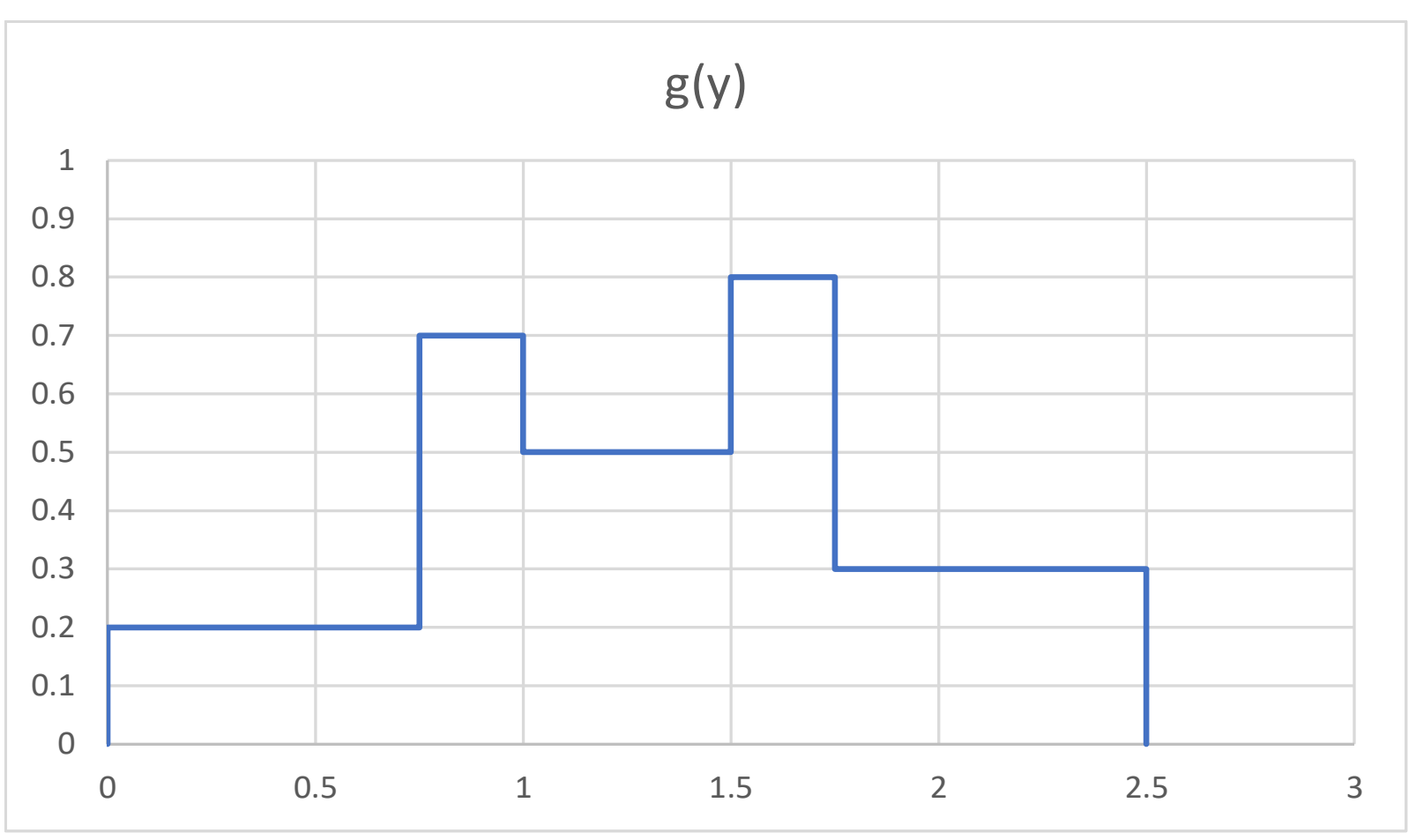

and the partition of the unit square implied by this $g()$ is

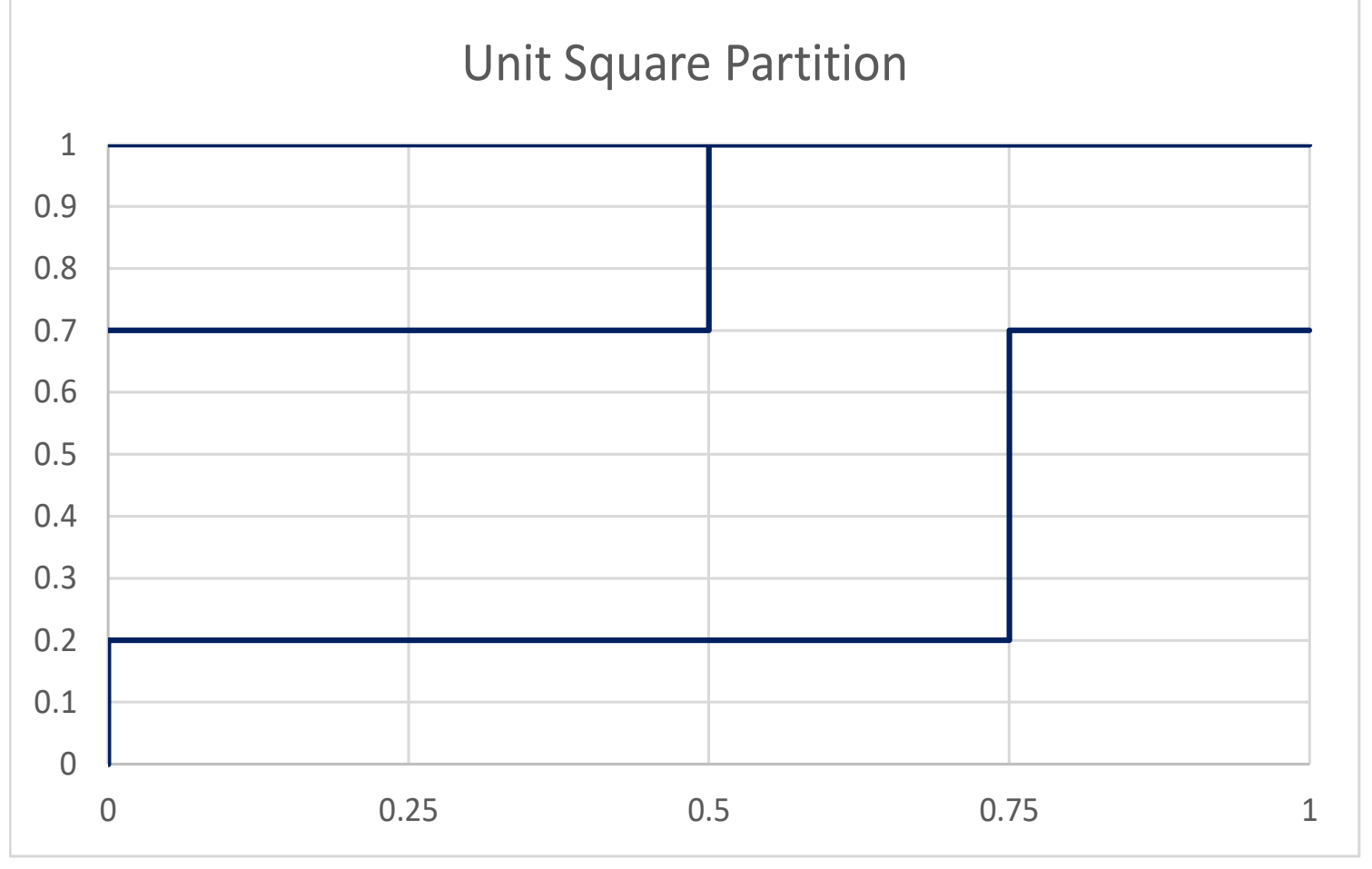

**(2) Cumulative distribution functions**. As noted above, the cdf of *every* random variable, whether discrete, continuous, or mixed, provides a valid seed function.

The cdf of any discrete valued random variable is a good example of a "multiple jump" seed function. For example, a "binomial" seed function would have $q_k = k \in \{0, 1, \ldots, n\}$ and $\alpha_k = \binom{n}{k} p^k (1-p)^{n-k}$ for some $n \in \mathbb{N}$ and $0 < p < 1$, and a "Poisson" seed function would have $q_k = k \in \mathbb{N}_0$ and $\alpha_k = e^{-\lambda}(\lambda^k / k!)$ for some $\lambda > 0$. The cdf of a standard normal used in Example 4.3 is an archetypal example of the use of the cdf of a continuous random variable as a seed function.

**(3) Other increasing functions**. Any function that *could* be used as a cdf is a valid seed function. Restricting attention just to continuous functions, these range from those that are smooth enough to be elements of $C^1$, e.g.,

$$H(y) \equiv \begin{cases} 0 & \text{if} \quad y < 0, \\ \sin^2\left(\frac{1}{2}\pi(\frac{y}{L})\right) & \text{if} \quad 0 \le y < L, \\ 1 & \text{if} \quad y \ge L, \end{cases}$$

for some $L > 0,$ to those that are continuous but not absolutely continuous, e.g., the Cantor function.

**(4) Functions that are unit interval increasing but not increasing**. I'll give two examples. The first example is artificial, but the second is of some importance.

**Example 1**. Let $F()$ be a continuous and strictly increasing cdf. Define the "switch" function $s()$ as

$$s(y) \equiv \lfloor y \rfloor + \lceil y \rceil - y.$$

For any integer $k$, as $y$ increases on the interval $(k, k+1)$, $s(y)$ *decreases* from $k+1$ to $k$. The function

$$H(y) \equiv F(s(y))$$

is unit interval increasing but not increasing. In fact, it's strictly decreasing on every interval $(k, k+1)$ where $k \in \mathbb{Z}$. When $y \in \mathbb{Z}$, $H()$ has a jump discontinuity of size $F(y+1) - F(y-1)$.

**Example 2**. Let $Y \sim g()$ where

$$g(y) = \frac{1 - \cos(2\pi y)}{2\pi^2 y^2}, \tag{6.1}$$

shown below.

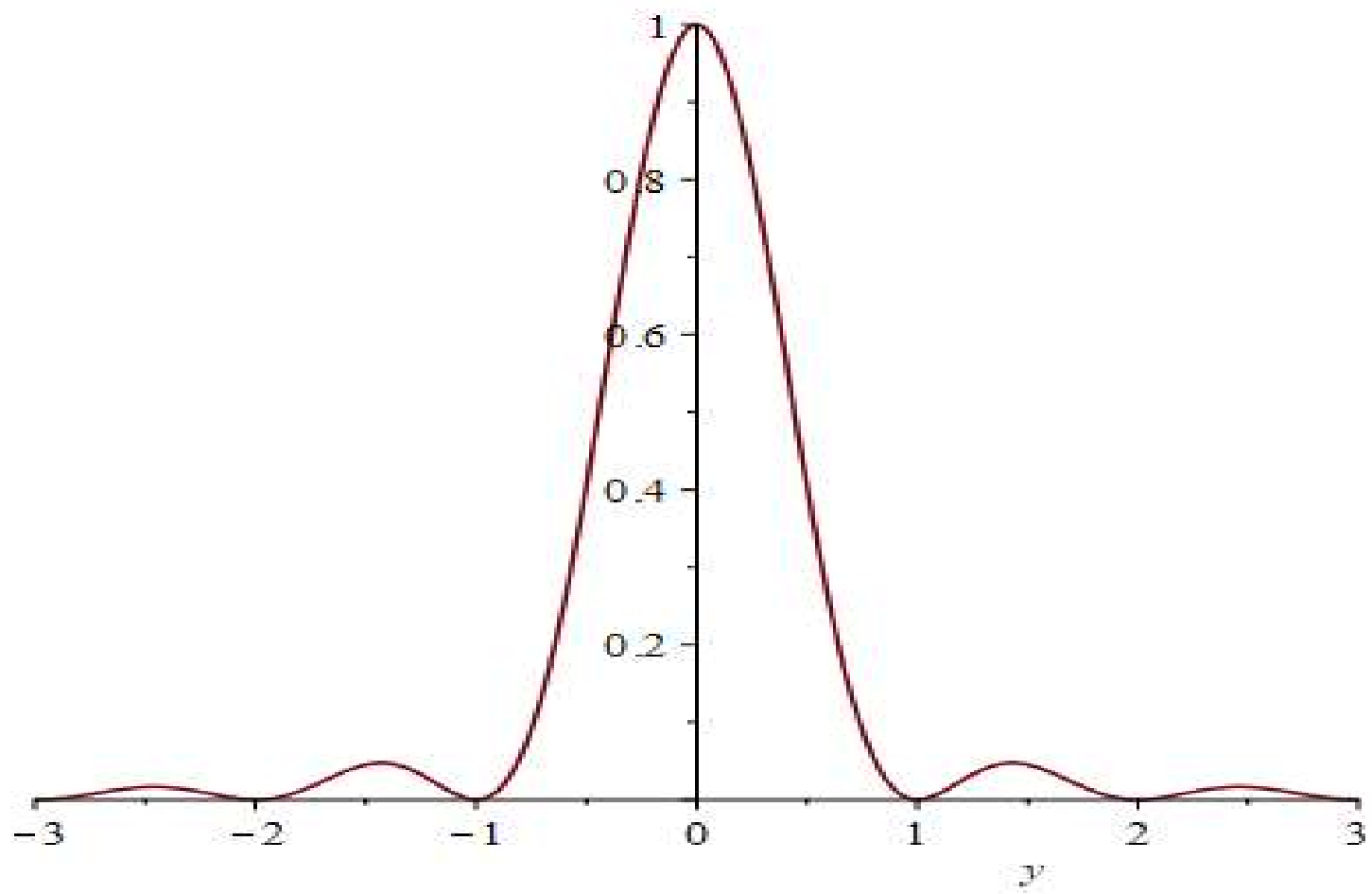

This is the pdf of $Y \equiv \log_b(X)$ where $X$ is Whittaker's random variable, discussed in [2]. Fourier analysis may be used to show that $Y$ is u.d. mod 1, and Maple assures me that the associated seed function is given by

$$H(y) = \frac{\Psi'(-y)\sin(\pi y)^2}{\pi^2} \tag{6.2}$$

where $\Psi()$ denotes the digamma function. This seed function is plotted below.

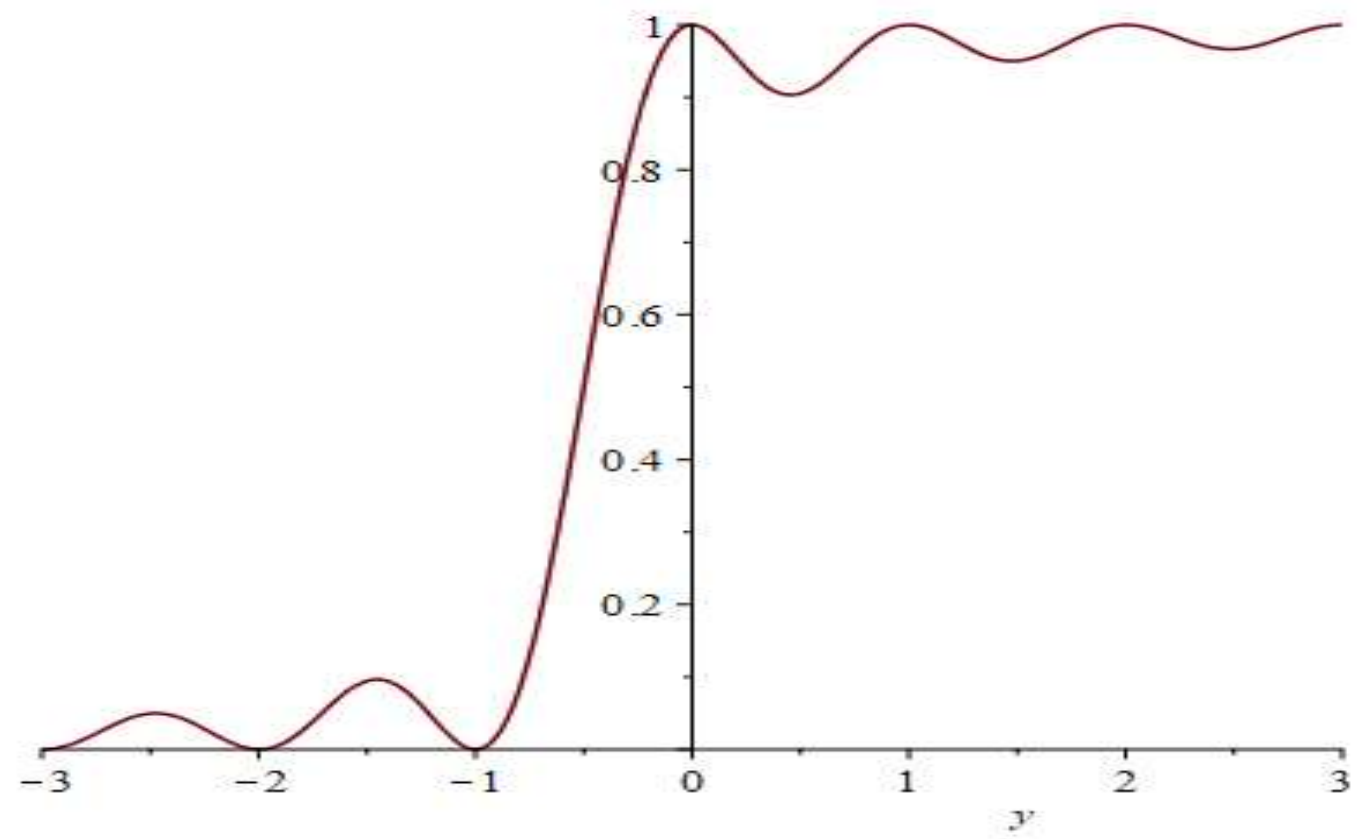

## 7. Piecewise Linear Seed Functions.

In this section we examine a simple family of seed functions that have a surprising property.

Let the parameter $r > 0$ be fixed and define $s \equiv r^{-1}$. Then a one parameter family of piecewise linear seed functions is defined as

$$H(y) = \begin{cases} 0 & \text{if} \quad y < 0, \\ ry & \text{if} \quad 0 \le y < s, \\ 1 & \text{if} \quad y \ge s. \end{cases} \tag{7.1}$$

The implied pdf $g()$ is triangular if $r = 1$ , and trapezoidal if $r \neq 1$. Specifically, if $r = 1$ then

$$g(y) = \begin{cases} y & \text{if} \quad 0 \le y < 1, \\ 2 - y & \text{if} \quad 1 \le y < 2, \end{cases} \tag{7.2}$$

and $g(y) = 0$ for other values of $y$. If $r < 1$ then

$$g(y) = \begin{cases} ry & \text{if} \quad 0 \le y < 1, \\ r & \text{if} \quad 1 \le y < s, \\ 1 + r - ry & \text{if} \quad s \le y < 1 + s, \end{cases} \tag{7.3}$$

and $g(y) = 0$ for other values of $y$. Finally, if $r > 1$ then

$$g(y) = \begin{cases} ry & \text{if} \quad 0 \le y < s, \\ 1 & \text{if} \quad s \le y < 1, \\ 1 + r - ry & \text{if} \quad 1 \le y < 1 + s, \end{cases} \tag{7.4}$$

and $g(y) = 0$ for other values of $y$. The later two cases are illustrated below.

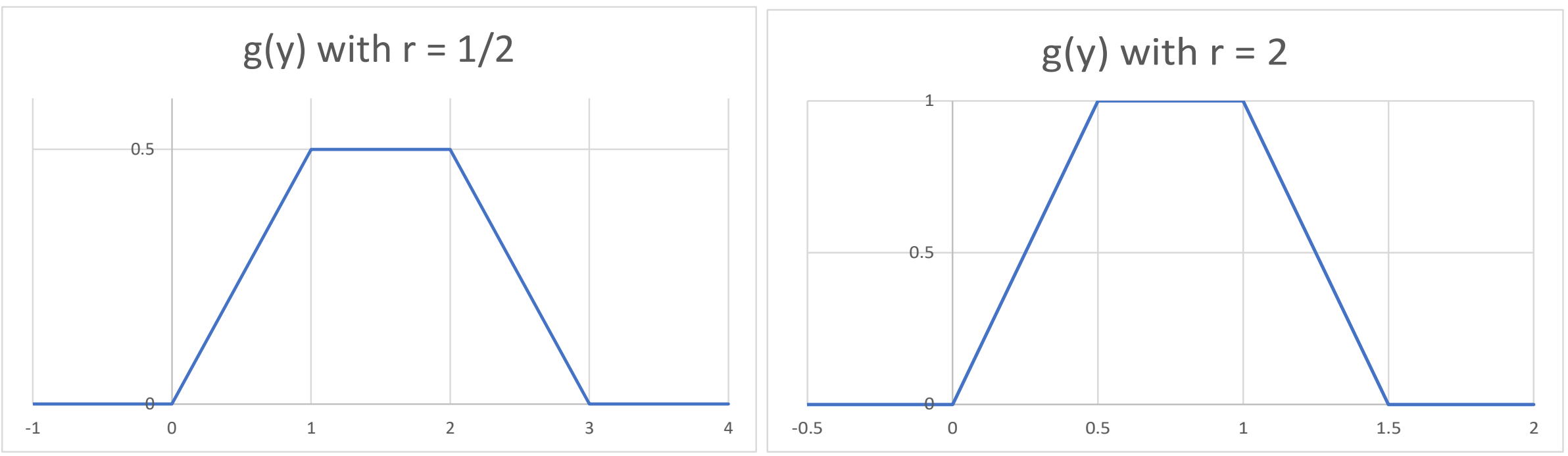

The implied pdf $f(x)$ in either of the two $r \neq 1$ cases has a three part formula on its support. Assuming a base $b$, when $r < 1$

$$\begin{aligned} 1 \le x < b \quad &\Rightarrow \quad f(x) = \frac{\Lambda_b}{x}(r \log_b(x)), \\ b \le x < b^s \quad &\Rightarrow \quad f(x) = \frac{\Lambda_b r}{x}, \\ b^s \le x < b^{(1+s)} \quad &\Rightarrow \quad f(x) = \frac{\Lambda_b}{x}(1 + r - r \log_b(x)), \end{aligned} \tag{7.5}$$

and $f(x) = 0$ elsewhere. The formula when $r > 1$ is similar:

$$
\begin{aligned}
1 \le x < b^s \quad &\Rightarrow \quad f(x) = \frac{\Lambda_b}{x}(r\log_b(x)), \\
b^s \le x < b \quad &\Rightarrow \quad f(x) = \frac{\Lambda_b}{x}, \\
b \le x < b^{(1+s)} \quad &\Rightarrow \quad f(x) = \frac{\Lambda_b}{x}(1 + r - r\log_b(x)),
\end{aligned}
\tag{7.6}
$$

and $f(x) = 0$ elsewhere. A graph of $f(x)$ when $r = 2$ and $b = 5$ is shown below:

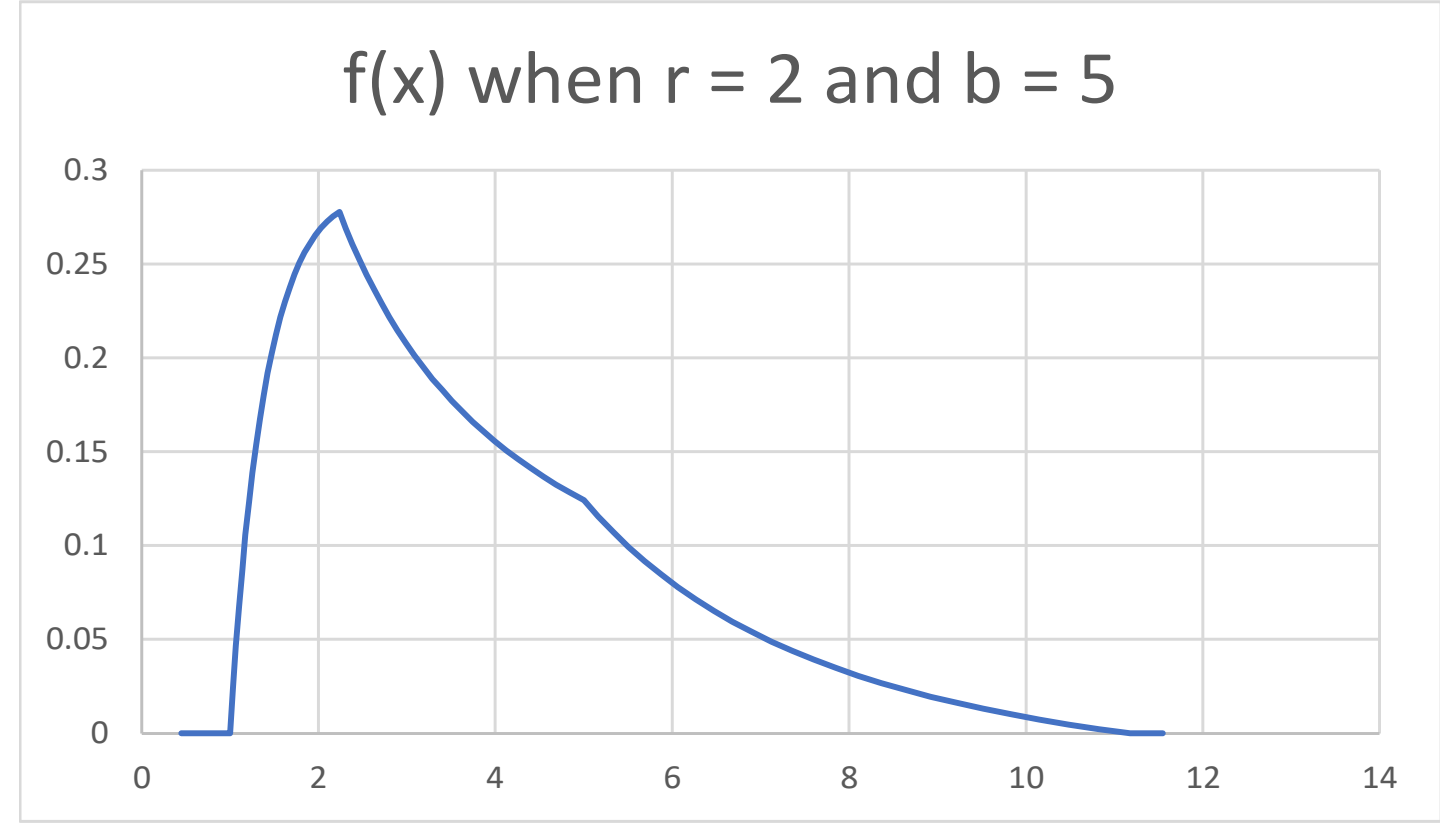

This graph has kinks at $x = 1,\ x = \sqrt{5},\ x = 5,$ and at $x = 5^{3/2} = 5\sqrt{5}$.

Suppose that $X \sim f()$, where $f()$ is given by either eq. (7.5) or (7.6). Then $X$ is $b$-Benford by construction, and $Y \equiv \log_b(X) \sim g()$, where $g()$ is given by either eq. (7.3) or (7.4). Let $a > 1$ be an alternative possible base, let $Y_a \equiv \log_a(X)$, and let $g_a()$ denote the pdf of $Y_a$. From eq. (5.3),

$$
g_a(y) = pg(py), \tag{7.7}
$$

where $p \equiv \ln(a)/\ln(b)$. As noted above, $Y_a$ is generally not u.d. mod 1, and hence $X$ is not generally $a$-Benford. However, if we choose $a$ cleverly, we may show that $X$ is $a$-Benford.

The key is to choose $a$ so that $p = r^{-1} = s$. That is,

$$
\frac{\ln(a)}{\ln(b)} = s \quad \Leftrightarrow \quad a = b^s. \tag{7.8}
$$

Note that $r > 1$ implies that $a < b$ and $r < 1$ implies that $a > b$.

Suppose that $r < 1$. Combining eqs. (7.3) and (7.7), we find that

$$
g_a(y) = sg(sy) = \begin{cases} srsy & \text{if} \quad 0 \le sy < 1, \\ sr & \text{if} \quad 1 \le sy < s, \\ s(1 + r - rsy) & \text{if} \quad s \le sy < 1 + s, \end{cases}
$$

and $g_a(y) = 0$ for other values of $y$. After cancellation, this expression settles down to

$$g_a(y) = \begin{cases} sy & \text{if} \quad 0 \le y < r, \\ 1 & \text{if} \quad r \le y < 1, \\ 1+s-sy & \text{if} \quad 1 \le y < 1+r, \end{cases} \tag{7.9}$$

and $g_a(y) = 0$ for other values of $y$.

Similarly, if $r > 1$ we find that

$$g_a(y) = \begin{cases} sy & \text{if} \quad 0 \le y < 1, \\ s & \text{if} \quad 1 \le y < r, \\ 1+s-sy & \text{if} \quad r \le y < 1+r, \end{cases} \tag{7.10}$$

and $g_a(y) = 0$ for other values of $y$.

Comparison of eqs. (7.9) and (7.10) with eqs. (7.3) and (7.4) reveals a remarkable symmetry: the two pairs of equations are identical except that the roles of $r$ and $s$ are reversed. In fact, it is easy to see that $g_a(y) = H_a(y) - H_a(y-1)$ where $H_a()$ is the piecewise linear seed function

$$H_a(y) \equiv \begin{cases} 0 & \text{if} \quad y < 0, \\ sy & \text{if} \quad 0 \le y < r, \\ 1 & \text{if} \quad y \ge r. \end{cases} \tag{7.11}$$

This is identical to the definition of $H()$ except that the roles of $r$ and $s$ are reversed.

In this discussion we took $r$ and $b$ as primitive and defined $a$ in terms of $r$ and $b$. But we could equally well have taken $a$ and $b$ as primitive and defined $r$ in terms of $a$ and $b$. To be specific, given $a > 1$, $b > 1$, and $a \ne b$, we define $r \equiv \ln(b)/\ln(a)$.

## 8. An Introduction to Product Random Variables.

The fact that the random variable $X$ defined in the previous section is both $a$-Benford and $b$-Benford (for any alternative bases $a$ and $b$ with $a \ne b$ that one may choose) is not so miraculous as it may appear. In fact, it's a consequence of the following proposition, which is a corollary of Berger and Hill's Theorem 8.12, page 188.

**Proposition 8.1**. Suppose that $X_a$ and $X_b$ are independent positive random variables. If $X_a$ is $a$-Benford and $X_b$ is $b$-Benford, then the product $X \equiv X_a X_b$ is both $a$-Benford and $b$-Benford.

With the random variables of Proposition 8.1, let $f_a()$ denote the pdf of $X_a$, $f_b()$ denote the pdf of $X_b$, and let $f()$ denote the pdf of the product random variable $X = X_a X_b$. Then

it may be shown that

$$f(x) = \int_{0+}^{\infty} \frac{1}{\xi} f_a(\xi) f_b\left(\frac{x}{\xi}\right) d\xi. \tag{8.1}$$

To take the simplest possible example of this formula, suppose

$$f_a(x) = \frac{\Lambda_a}{x} \quad \text{if} \quad 1 \le x \le a, \qquad f_a(x) = 0 \quad \text{elsewhere,} \tag{8.2a}$$

and

$$f_b(x) = \frac{\Lambda_b}{x} \quad \text{if} \quad 1 \le x \le b, \qquad f_b(x) = 0 \quad \text{elsewhere.} \tag{8.2b}$$

Note that $f_a(\xi) f_b(\eta) > 0$ iff $(\xi, \eta) \in [1, a] \times [1, b]$, and hence the support of $X$ is $[1, ab]$. For any $x \in [1, ab]$, note that

$$f_a(\xi) f_b\left(\frac{x}{\xi}\right) = \frac{\Lambda_a}{\xi} \cdot \frac{\Lambda_b}{x/\xi} = \frac{\Lambda_a \Lambda_b}{x} \tag{8.3}$$

for any $\xi$ where $f_a(\xi)$ and $f_b(x/\xi)$ are both positive. Let $\Lambda \equiv \Lambda_a \Lambda_b$. Without loss of generality, assume that $a < b$. The limits of the integral of eq. (8.1) depend on which of the disjoint ranges $(1, a)$, $(a, b)$, and $(b, ab)$ contains $x$, so the pdf $f()$ has a three part formula. Leaving the details to the reader, we may show the following: for $x \in (1, a)$,

$$f(x) = \int_1^x \frac{\Lambda}{x} \cdot \frac{1}{\xi}\, d\xi = \frac{\Lambda}{x} \ln(x), \tag{8.4a}$$

for $x \in (a, b)$,

$$f(x) = \int_1^a \frac{\Lambda}{x} \cdot \frac{1}{\xi}\, d\xi = \frac{\Lambda}{x} \ln(a), \tag{8.4b}$$

and finally, for $x \in (b, ab)$,

$$f(x) = \int_{x/b}^a \frac{\Lambda}{x} \cdot \frac{1}{\xi}\, d\xi = \frac{\Lambda}{x}[\ln(a) - \ln(x/b)] = \frac{\Lambda}{x} \ln\left(\frac{ab}{x}\right). \tag{8.4c}$$

The assumption that $a < b$ implies that $r \equiv \ln(b)/\ln(a) > 1$. The reader may confirm that eqs. (8.4) are identical to eq. (7.6).

Let's generalize some of these results. Suppose that $X_a \sim f_a()$ is $a$-Benford, that $X_b \sim f_b()$ is $b$-Benford, and that $X_a$ and $X_b$ are independent. As $X_a$ is $a$-Benford, there exists a seed function $H_a()$, a pdf $g_a()$, and a u.d. mod 1 random variable $Y_a$ associated with $(X_a, f_a())$ as shown in the following schematic:

$$H_a() \leftrightarrow (Y_a, g_a()) \overset{a}{\leftrightarrow} (X_a, f_a()). \tag{8.5}$$

Similarly, there exists a seed function $H_b()$, a pdf $g_b()$, and a u.d. mod 1 random variable $Y_b$

associated with $(X_b, f_b())$ as shown in the schematic

$$H_b() \leftrightarrow (Y_b, g_b()) \overset{b}{\leftrightarrow} (X_b, f_b()). \tag{8.6}$$

Under these assumptions, the product $X \equiv X_a X_b$ is both $a$-Benford and $b$-Benford. Let $f()$ denote $X$'s pdf, given by eq. (8.1). As $X$ is both $a$-Benford and $b$-Benford, there exist seed functions $H_a^*()$ and $H_b^*()$ associated with $(X, f())$ according to the following schematics:

$$\begin{aligned} H_a^*() &\leftrightarrow (Y_a^*, g_a^*()) \overset{a}{\leftrightarrow} (X, f) \\ H_b^*() &\leftrightarrow (Y_b^*, g_b^*()) \overset{b}{\leftrightarrow} (X, f). \end{aligned} \tag{8.7}$$

The question then arises: what is the relationship between $(H_a(), H_b())$ and $(H_a^*(), H_b^*())$? So far, we have a single example. If

$$H_a(y) = H_b(y) = \begin{cases} 0 & \text{if} \quad y < 0, \\ 1 & \text{if} \quad y \geq 0, \end{cases} \tag{8.8}$$

then

$$H_a^*(y) = \begin{cases} 0 & \text{if} \quad y < 0, \\ sy & \text{if} \quad 0 \leq y < r, \\ 1 & \text{if} \quad y \geq r, \end{cases} \quad \text{and}$$

$$H_b^*(y) = \begin{cases} 0 & \text{if} \quad y < 0, \\ ry & \text{if} \quad 0 \leq y < s, \\ 1 & \text{if} \quad y \geq s, \end{cases}$$

where $s \equiv \ln(a)/\ln(b)$ and $r \equiv s^{-1}$.

Here's a related question. Let $c > 1$, define $Y_c^* \equiv \log_c(X)$, and let $g_c^*()$ denote the pdf of $Y_c^*$. We know that $Y_c^*$ is not generally u.d. mod 1, although it is when $c = a$ or when $c = b$. It is of some interest to examine the dependence of $g_c^*()$ on $c$. As $X = X_a X_b$, it follows that

$$\begin{aligned} Y_c^* &= \log_c(X_a) + \log_c(X_b) \\ &= \frac{\ln(a)}{\ln(c)} \frac{\ln(X_a)}{\ln(a)} + \frac{\ln(b)}{\ln(c)} \frac{\ln(X_b)}{\ln(b)} = \rho(a, c) Y_a + \rho(b, c) Y_b \end{aligned} \tag{8.9}$$

where $\rho(d, c) \equiv \ln(d)/\ln(c)$ for any $d > 1$. If $c = a$, eq. (8.9) becomes

$$Y_a^* = Y_a + \rho(b, a) Y_b, \tag{8.10}$$

and if $c = b$ it becomes

$$Y_b^* = \rho(a, b) Y_a + Y_b. \tag{8.11}$$

(In our earlier notation, $\rho(a, b) = s$ and $\rho(b, a) = r$.)

Equation (8.9) implies that $g_c^*()$ is the convolution of the pdfs of $\rho(a, c) Y_a$ and $\rho(b, c) Y_b$. For example, suppose that $c = a$ and (to simplify notation) we rewrite eq. (8.10) as

$$Y_a^* = Y_a + \rho Y_b \qquad \text{where} \qquad \rho \equiv \rho(b, a).$$

The pdf of $\rho Y_b$ is $\rho^{-1} g_b(\rho^{-1} y)$, and it follows that

$$g_a^*(y) = \rho^{-1} \int_{-\infty}^{\infty} g_a(\xi) g_b\big(\rho^{-1}(y - \xi)\big) d\xi. \tag{8.12}$$

Continuing this example, suppose that both $Y_a$ and $Y_b$ are uniformly distributed on $[0, 1]$. Then the support of $(Y_a, \rho Y_b)$ is $[0, 1] \times [0, \rho]$. This implies that our formulas for $g_a^*(y)$ will have a three part form. There are two cases to consider. If $\rho > 1$,

$$g_a^*(y) = \begin{cases} \rho^{-1} y & \text{if} \quad 0 \le y < 1, \\ \rho^{-1} & \text{if} \quad 1 \le y < \rho, \\ 1 + \rho^{-1} - \rho^{-1} y & \text{if} \quad \rho \le y < 1 + \rho, \end{cases} \tag{8.13}$$

and $g_a^*(y) = 0$ elsewhere. If $\rho < 1$,

$$g_a^*(y) = \begin{cases} \rho^{-1} y & \text{if} \quad 0 \le y < \rho, \\ 1 & \text{if} \quad \rho \le y < 1, \\ 1 + \rho^{-1} - \rho^{-1} y & \text{if} \quad 1 \le y < 1 + \rho, \end{cases} \tag{8.14}$$

and $g_a^*(y) = 0$ elsewhere. We see that eqs. (8.13) and (8.14) are just a rederivation of eqs. (7.9) and (7.10).

These equations give us a start on a "Benford analysis" of the product random variable $X$, but further progress requires some results from Fourier analysis and must be deferred to another paper.

## 9. Concluding Remarks.

One of my aims when I began this research was to find an easy way to generate lots of examples of pdfs of Benford random variables. With regard to this goal, at least, I can say that I succeeded beyond my wildest dreams. The method of seed functions allow one to generate an unlimited number of examples of pdfs of Benford random variables.

One use for these tools is in an investigation into some issues of base dependence in Benford random variables. Suppose $X$ is $b$-Benford. We know that $X$ is generally not $c$-Benford when $c \neq b$. That is, the random variable $Y_c \equiv \log_c(X)$ is not generally u.d. mod 1, and the pdf of $\langle Y_c \rangle$, say $\widetilde{g}_c(u)$, is not generally identically equal to 1 for almost all $u \in [0, 1)$. In a companion paper [1], I combined the tools outlined here with Fourier analysis to examine the dependence of $\widetilde{g}_c()$ on $c$ for $X$ in several families of Benford random variables. I also used these tools to examine the dependence of $\widetilde{g}_c()$ on $c$ for the product random variable $X \equiv X_a X_b$, where $X_a$ is $a$-Benford and $X_b$ is $b$-Benford.

The reader is invited to find his/her own uses for these tools, either in research or in teaching.

## References

[1] Benford, Frank. "Base Dependence of Benford Random Variables." arXiv:1702.01644, 2017.

[2] Benford, Frank. "Fourier Analysis and Benford Random Variables." arXiv:2006.07136, 2020.

[3] Berger, Arno, and Theodore Hill. *An Introduction to Benford's Law*. Princeton University Press, Princeton and Oxford, 2015.